\documentclass[12pt]{amsart}
\usepackage{xspace}
\usepackage[curve,cmtip,matrix,arrow]{xy}
\UseTips 

\usepackage[bookmarks=true]{hyperref}

\newtheorem{theorem}{Theorem}[section]
\newtheorem{corollary}[theorem]{Corollary}
\newtheorem{lemma}[theorem]{Lemma}
\newtheorem{proposition}[theorem]{Proposition}

\newtheorem*{PL}{Poincare Lemma}
\theoremstyle{definition}
\newtheorem{definition}[theorem]{Definition}
\theoremstyle{remark}
\newtheorem{remark}[theorem]{Remark}

\numberwithin{equation}{section}

\newcommand{\QQ}{\mathbb{Q}}
\newcommand{\ZZ}{\mathbb{Z}}
\newcommand{\FF}{\mathbb{F}}

\renewcommand{\to}{\xrightarrow}
\newcommand{\To}{\longrightarrow}
\newcommand{\cal}[1]{\mathcal{#1}}
\newcommand{\HHom}{\mathcal{H}om}

\newcommand{\HH}{\cal{H}}
\newcommand{\Hyper}{\mathbb{H}}
\newcommand{\OO}{\cal{O}}
\newcommand{\Spec}{\operatorname{Spec}}
\newcommand{\codim}{\operatorname{codim}}
\newcommand{\Hom}{\operatorname{Hom}}

\newcommand{\wt}{\widetilde}
\newcommand{\Omtilde}{\wt{\Omega}}

\newcommand{\deRham}{\mbox{de$\,$Rham}\xspace}

\begin{document}

\title{Cartier Isomorphism for Toric  Varieties}
\author{Manuel Blickle}
\address{Univ. of Michigan, 525 East University Avenue, Ann Arbor, MI-48109}
\email{blickle@math.lsa.umich.edu}

\thanks{To appear in Journal of Algebra}

\date{June 11, 2000}




\maketitle


\section{Introduction}

In \cite{Cartier57}, Cartier introduces an operation on the \deRham complex for
varieties defined over a field of positive characteristic. In the case of a
smooth variety this yields a complete description of the cohomology of the
\deRham complex. To formulate the result let $X$ be a smooth variety over a
perfect field then for all $a \geq 0$ the Cartier operator
\[{C: \HH^a(F_*\Omega^\bullet_X) \to{}  \Omega^a_X}  \]
is an isomorphism. Here $F: X \to{} X$ denotes the Frobenius morphism on $X$
and $\HH^a$ denotes the $a^{\text{th}}$ cohomology sheaf of
$F_*\Omega^{\bullet}_X$. If the variety is not smooth, not much is known about
the properties of the Cartier operator and the poor behaviour of the \deRham
complex in this case makes its study difficult. If one substitutes the \deRham
complex with the Zariski-\deRham complex the situation is better. For example,
the Zariski differentials, though not locally free, are reflexive and there is
a natural duality pairing between them.

We show how to extend the Cartier operator in a natural way to the Zariski
differentials. Using a description of the Zariski-\deRham complex due to
Danilov \cite{Danilov78} we show that this newly defined Cartier operator is an
isomorphism for toric varieties. Moreover, it is induced by a split injection
\[
\xymatrix@1{ 0 \ar[r] &{\Omtilde^a_X} \ar[r] &{F_*\Omtilde^a_X.}}
\]
As described in \cite{Anders} such a result yields the Bott vanishing theorem
and the degeneration of the Hodge to \deRham spectral sequence for projective
toric varieties.

Finally we will give an obstruction to the surjectivity of the Cartier map to
show that in general it can not be an isomorphism.

The author would like to thank Karen Smith for invaluable discussions
and the referee for a careful reading and valuable comments.

\subsection{Notation and Generalities}

Fix a perfect field $k$ of characteristic $p>0$. All schemes will be noetherian
$k$-schemes unless otherwise stated. $X$ denotes a noetherian $k$-scheme. By an
$\OO_X$-module we always mean a coherent sheaf of $\OO_X$-modules.

\begin{definition}
  The \emph{absolute Frobenius} of $X$ is the endomorphism
  \[F_X :X \To X\]
  which is the identity on the level of the topological space and the
  $p^{\text{th}}$ power map on the structure sheaf.
\end{definition}

For $x \in \OO _X$ and $\lambda \in k$ one has $F^*_X(\lambda
x)=\lambda^px^p=F_k(\lambda)F^*_X(x)$ where $F_k$ is just the Frobenius
morphism for $k$. Thus $F_X$ is \emph{not} a morphism of $k$-schemes unless
$k=\FF_p$. Setting $X'=X\times_{F_k} k$  we get by the universal property of
the fiber product that $F_X$ can be factored through $X'$.
\begin{equation*}
\xymatrix@=9pt{ {X} \ar[rrd]^{F}\ar@/_.5pc/[rrddd] \ar@/^1pc/[rrrrd]^{F_X} \\
                                 && {X'} \ar[rr]^{\pi_1}\ar[dd] && {X} \ar[dd] \\ \\
                                 && {k} \ar[rr]^{F_k}           && {k} }
\end{equation*}
This map is called the relative Frobenius and we simply denote it by $F: X
\to{} X'$. Note that as $\FF_p$-schemes $X$ and $X'$ are isomorphic. For an
$\OO_X$-module $\cal{M}$ we denote the pushforward under $F$ by $F_*\cal{M}$.

For every map $\phi: X \to{} Y$ of noetherian $k$-schemes we get a commutative
diagram:
\[
\xymatrix{ {X} \ar^F[r] \ar_{\phi}[d]  &{X'} \ar^{\phi'}[d] \\
           {Y} \ar^F[r]                &{Y'}
         }
\]
This follows by construction of the relative Frobenius as the unique solution
to a mapping problem and the fact that we have a square like this for the
absolute Frobenius. Thus, in particular, the functor $F_*$ commutes with
$\phi_*$.

For a $k$-algebra $R$ we adopt the same notation: The (absolute) Frobenius
$F_R: R \to{} R$ is just the $p^{\text{th}}$ power map. Setting $R'=R \otimes_F
k$ we get similarly a relative ($k$-linear) Frobenius map $F: R' \to{} R$.
Again one sees that $R$ and $R'$ are isomorphic as rings, only their
$k$-algebra structure is different. By $F_*M$ we denote the $R$-module $M$
viewed as an $R'$-module via the relative Frobenius.

\subsection{The smooth case}

In this section we will review the construction in the smooth case which we
will extend later to the Zariski-\deRham complex on normal varieties.

The object of study is the algebraic \deRham complex
$\Omega^{\bullet}_{X/k}=\Omega^\bullet_X$ (cf. \cite{EV}). The key observation
is that the differential on $F_*\Omega^\bullet_X$ is $\OO_{X'}$-linear by the
Leibniz rule: $dF^*(x\otimes1)=dx^p=0$ for all $x\in \OO_X$. This gives the
cohomology objects of the de Rham complex $\HH^a(F_*\Omega^\bullet_{X/S})$ the
structure of $\OO_{X'}$-modules.

\begin{proposition}[Cartier, \cite{Cartier57}] Let $X$ be a noetherian $k$-scheme.
  There is a unique map of $\OO_{X'}$-modules
  \[ C^{-1}: \Omega^1_{X'}\To \HH^1(F_*\Omega^\bullet_X)\]
  called the \emph{inverse Cartier operator} such that for a local
  section $x\in \OO_X$ one has
  \[
  C^{-1}(d(x\otimes 1)) = x^{p-1}dx \text{ \ in } \HH^1.
  \]
  For all $a \geq 0$ this $C^{-1}$ induces maps
  \[
  \bigwedge^a C^{-1}: \Omega^a_{X'} \To \HH^a(F_*\Omega^\bullet_X).
  \]
\end{proposition}

The following result completely describes the cohomology of the \deRham complex
on a smooth $k$-scheme.
\begin{theorem}[Cartier, \cite{Cartier57}, see also
\cite{Katz}]\label{smooth:cartier} Let $X$ be a smooth $k$-scheme. Then the
 inverse Cartier operator is an isomorphism for all $a \geq 0$. I.e.
 \[
 \HH^a(F_*\Omega^\bullet_X) \cong \Omega^a_{X'}.
 \]
\end{theorem}

\begin{remark}\label{Cart:DI}
The most striking application of the Cartier isomorphism is as the main
ingredient of Deligne and Illusie's algebraic proof of the Kodaira vanishing
theorem. For $X$ smooth and under an additional assumption ($X$ lifts to
characteristic $p^2$ and $p \geq \dim X$) they show in \cite{DI} that the
inverse Cartier operator is induced by a map
\[ \bigoplus \Omega_{X'}^a[-a] \to{} F_*\Omega^\bullet_X \]
in the derived category of $\OO_{X'}$-modules. Taking hypercohomology one sees
that this implies the degeneration of the Hodge to \deRham spectral sequence
$H^a(\Omega_X^b) \Rightarrow \Hyper^{a+b}(\Omega^\bullet_X)$. Such a result
implies the Kodaira vanishing theorem as shown by Esnault and Viehweg in
\cite{EVLog}. For details on this circle of ideas we refer the reader to the
original article \cite{DI} and the excellent book \cite{EV}.
\end{remark}

\subsection{Reflexive sheaves}

In our approach to generalize the Cartier operator reflexive sheaves play a
central role. Here we review some properties and for completeness most of the
proofs are given, too.

 Let $X$ be a normal scheme and $\cal{M}$ be a
$\OO_X$-module. We denote the dual of $\cal{M}$ by $\cal{M}^* =
\HHom_{\OO_X}(\cal{M},\OO_X)$. There is a natural map from $\cal{M}$ to its
double dual $\cal{M}^{**}$ which sends a section $s \in \cal{M}$ to the map
$f_s$ whose value on the map $\lambda: \cal{M} \to{} \OO_X$ is
$f_s(\lambda)=\lambda(s)$. The module $\cal{M}$ is called reflexive if the
natural map to its double dual is an isomorphism. The module $\cal{M}^{**}$ is
also called the reflexive hull of $\cal{M}$. It is reflexive \cite{EvGr} and
universal with respect to this property: every map $\cal{M} \to{} \cal{N}$ with
$\cal{N}$ reflexive factors canonically through the natural map from $\cal{M}$
to $\cal{M}^{**}$.

Recall the following local version of the adjointness of $f^*$ and $f_*$ for a
map of schemes $f: X \to{} Y$.

\begin{lemma}\label{localadj}
  Let $f: X \to{} Y$ be a map of schemes, $\cal{F}$ an $\OO_X$-module and $\cal{G}$
  an $\OO_Y$-module. There is a natural isomorphism
  \[ \Phi: f_* \HHom_{\OO_X}(f^*\cal{G},\cal{F}) \to{\ \cong\ } \HHom_{\OO_Y}(\cal{G},f_*\cal{F}).\]
\end{lemma}

\noindent To see there is a globally defined map let $\phi_U$ for $U\subseteq
Y$ be a local section of the left hand side. Composing the natural map
$\cal{G}|_U \to{} f_*f^*\cal{G}|_U=(f_*f^*\cal{G})|_U$ with
$f_*(\phi_U):f_*f^*\cal{G}|_U \to{} f_*\cal{F}|_U$ we get the desired local
section of the right hand side. Since no choices were made this defines the
global map $\Phi$. An inverse of $\Phi$ can be established similarly using the
natural map $f^*f_* \cal{F} \to{} \cal{F}$.

\begin{lemma}\label{ReflExt} Let $X$ be a normal scheme and let
$i:U \hookrightarrow X$ the inclusion of an open set $U$ such that
$\codim(X,X-U) \geq 2$. If $\cal{M}$ is a reflexive $\OO_U$-module then
$i_*\cal{M}$ is the unique reflexive $\OO_X$-module which agrees with $\cal{M}$
on $U$.
\end{lemma}
\begin{proof}
First we show that $i_*\cal{M}$ is reflexive by showing that $i_*$ commutes
with dualizing. For this we need the easily verifiable facts that for the given
setting $i_* \OO_U = \OO_X$ and $i^*i_*\cal{M} = (i_* \cal{M})|_U = \cal{M}$.
Now use the local adjointness Lemma \ref{localadj} in the following
calculation.
\begin{equation*}
  \begin{split}
  (i_*\cal{M})^* &= \HHom_{\OO_X}(i_*\cal{M},\OO_X) =
  \HHom_{\OO_X}(i_*\cal{M},i_*\OO_U)\\
               &= i_* \HHom_{\OO_U}(i^*i_*\cal{M},\OO_U) \\
               &= i_* \cal{M}^*
  \end{split}
\end{equation*}
Thus $(i_* \cal{M})^{**} \cong i_*\cal{M}^{**} \cong i_*\cal{M}$ which shows
that $i_* \cal{M}$ is reflexive. \\ For the uniqueness part let $\cal{N}$ be a
reflexive $\OO_X$-module such that its restriction to $U$ is $\cal{M}$. A
similar application of the adjointness lemma shows that
\begin{equation*}
  \begin{split}
    \cal{N}^* &= \HHom_{\OO_X}(\cal{N},\OO_X)
               = i_*\HHom_{\OO_U}(i^* \cal{N},\OO_U)\\
              &= i_* (\cal{N}|_U)^* = i_*\cal{M}^*
  \end{split}
\end{equation*}
and dualizing this equality we get $\cal{N}^{**} \cong i_* \cal{M}^{**}$ which
shows that $\cal{N} \cong i_* \cal{M}$ since both are reflexive.
\end{proof}

\begin{lemma}\label{refl:flat}
  Let $f:X \to{} Y$ be a flat morphism of noetherian normal schemes
  and let $\cal{M}$ be a reflexive $\OO_Y$-module. Then $f^* \cal{M}$ is a reflexive $\OO_X$-module.
\end{lemma}
\begin{proof} First we establish a natural map
$f^*\cal{M}^* \to{} (f^*\cal{M})^*$. Then we show by a local calculation that
this is an isomorphism whenever $f$ is flat. Thus dualizing commutes with $f^*$
and therefore $f^*$ preserves reflexivity.

For the first part notice that quite generally, for any, not necessarily flat,
$f$ and $\OO_Y$-modules $\cal{M}$ and $\cal{N}$ we get:
\[
f^*\HHom_{\OO_Y}(\cal{M},\cal{N}) \to{} f^*\HHom_{\OO_Y}(\cal{M},f_*f^*\cal{N})
\to{\ \cong\ } f^*f_*\HHom_{\OO_X}(f^*\cal{M},f^*\cal{N})
\]
The first map is composition with the natural map $\cal{N} \to{}
f_*f^*\cal{N}$ and the second map is the local adjointness Lemma
\ref{localadj}. Composing this with the natural map
$f^*f_*\HHom_{\OO_X}(f^*\cal{M},f^*\cal{N}) \to{}
\HHom_{\OO_X}(f^*\cal{M},f^*\cal{N})$ gives the desired map.

To see that this is an isomorphism when $f$ if flat we reduce to checking on
stalks. Thus it remains to show that for a flat map of local noetherian rings
$R \to{} S$ and finitely generated $R$-modules $M$ and $N$
\begin{equation}
 S \otimes_R \Hom_R(M,N) \to{\ \cong\ } \Hom_S(S \otimes_R M,S \otimes_R N). \tag{$*$} \label{local:eq}
\end{equation}
which sends $s \otimes f$ to the map $s'\otimes m \mapsto s's \otimes f(m)$ is
an isomorphism (this, of course, is the same as the map above). The case $M=R$
is easily verified since both sides are canonically isomorphic to $S \otimes_R
N$ and the above map is then the identity. Since tensor and $\Hom$ both commute
with finite direct sums (\ref{local:eq}) is an isomorphism for f.g. free
$R$-modules. Now let $F \to{} G \to{} M \to{} 0$ be a free presentation of $M$.
Applying $\Hom$ and tensor (exact!) to this presentation in both orders we get
a diagram:
\[
\xymatrix@=15pt{ 0 \ar[r] &{\Hom_R(M,N)\otimes S} \ar[r] \ar[d]
&{\Hom_R(G,N)\otimes S} \ar[r] \ar[d] &{\Hom_R(F,N)\otimes S} \ar[d] \\ 0
\ar[r] &{\Hom_S(M \otimes S,N \otimes S)} \ar[r] &{\Hom_S(G\otimes S,N \otimes
S)}\ar[r] &{\Hom_S(F\otimes S,N \otimes S)} }
\]
The vertical maps to the right are isomorphisms since $F$ and $G$ are free.
Thus by the 5-Lemma the left map is an isomorphism, too.
\end{proof}

\begin{corollary}\label{refl:fflat}
  Let $f:X \to{} Y$ be faithfully flat and $\cal{M}$ be an $\OO_X$-module. Then
  $\cal{M}$ is reflexive if and only if $f^*\cal{M}$ is reflexive.
\end{corollary}

\section{Cartier operator for normal schemes}

For schemes that are not smooth, the \deRham complex is not very well behaved
and it seems convenient to consider the Zariski-\deRham complex instead. In
order to define the sheaves of Zariski differentials we assume that the variety
$X$ is normal.

\subsection{Zariski-\deRham complex}

\begin{definition}
  Let $X$ be a normal $k$-scheme and let \[ i: U \hookrightarrow X \]
  be the inclusion of the smooth locus. The \emph{Zariski-\deRham complex} is the
  pushforward of the \deRham complex on $U$ and is denoted by
  $\Omtilde^\bullet_X=i_*\Omega^\bullet_U$.
\end{definition}

The objects of the Zariski-\deRham complex are just
$\Omtilde^a_X=i_*\Omega^a_U$ and are called the Zariski sheaves of differential
$a$-forms. Normality of $X$ guarantees that the Zariski sheaves are reflexive
by Lemma \ref{ReflExt}. Equivalently, one can define $\Omtilde^a_X$ as the
unique reflexive $\OO_X$-module which agrees with $\Omega^a_X$ on the smooth
locus (i.e $\Omtilde^a_X$ is the double dual or reflexive hull of $\Omega^a_X$
and the natural map $\Omega^a_X \to{} i_*i^*\Omega^a_X = \Omtilde^a_X$ is
identified with the double dualizing map). With this in mind one can take any
open smooth $U$ with codimension $X-U$ greater than or equal to two in the
above definition. The differential on $\Omtilde^\bullet_X$ is just the
pushforward $i_*(d)$ of the ordinary differential on $\Omega^{\bullet}_U$. The
Zariski sheaf of 1-forms $\Omtilde^1_X$ satisfies a universal property similar
to the universal property for the sheaf of K\"ahler differentials.

\begin{proposition}
  Every derivation $\delta:\OO_X \to{} \cal{M}$ to a reflexive $\OO_X$-module
  $\cal{M}$ factors uniquely through $i_*d:\OO_X \to{} \Omtilde^1_X$.
\end{proposition}
\begin{proof}
Observe that the derivation $i_*d$ is the compostition of the universal
derivation $d:\OO_X \to{} \Omega^1_X$ with the double dualizing map $\Omega^1_X
\to{} \Omtilde^1_X$. By the universal properties of these maps the proposition
is immediate.
\end{proof}

Note that, if $X=\Spec(R)$ is affine, the Zariski sheaf is the sheaf associated
to the finitely generated $R$-module
$\Gamma(X,\Omtilde^a_X)=\Gamma(U,\Omega^a_U)$. The Zariski-\deRham complex
shares many of the good properties of the \deRham complex on a smooth variety.
For example, if $X$ is $n$-dimensional, the next proposition shows that the
perfect pairing between $\Omega^a_U$ and $\Omega^{n-a}_U$ translates into a
perfect pairing for the Zariski sheaves.

\begin{proposition}[Danilov \cite{Danilov78}]
  Let $X$ be a n-dimensional normal $k$-scheme. There is a natural isomorphism
  $\Omtilde^a_X \rightarrow \HHom_{\OO_X}(\Omtilde^{n-a}_X,\Omtilde^n_X)$.
\end{proposition}
\begin{proof} Applying $i_*$ to the duality on the smooth locus and using the
local version of the adjointness of $i_*$ and $i^*$ we get
\[
\begin{split}
\Omtilde^a_X &= i_* \Omega^a_U = i_* \HHom_{\OO_U}(\Omega^{n-a}_U, \Omega^n_U) = i_* \HHom_{\OO_U}(i^*i_*\Omega^{n-a}_U, \Omega^n_U) \\
             &=\HHom_{\OO_X}(i_*\Omega^{n-a}_U,i_*\Omega^n_U)=\HHom_{\OO_X}(\Omtilde^{n-a}_X, \Omtilde^n_X)
\end{split}
\]
\renewcommand{\qed}{}
\end{proof}

\subsubsection{Stability under \'etale morphisms}

A morphism of sheaves $f: Y \to{} X$ is called \emph{\'etale} if it is
flat and unramified. Just as the sheaves of K\"{a}hler
differentials are preserved under \'{e}tale maps so are the Zariski sheaves.

\begin{lemma}\label{Zariskietale}
  Let $X$ be normal and $f: Y \to{} X$ \'{e}tale. Then
  \[ f^*\Omtilde^a_X \cong \Omtilde^a_Y. \]
\end{lemma}
\begin{proof}
Note that by \cite[Prop. 3.17(b)]{Milne} $Y$ is normal thus the Zariski sheaves are
defined on $Y$. Now
$f^*\Omtilde^a_X=f^*({\Omega^a_X}^{**})=(f^*\Omega^a_X)^{**}=\Omtilde^a_Y$
using Lemma \ref{refl:flat} and $f^*\Omega^a_X \cong \Omega^a_Y$ by \'etaleness
of $f$.
\end{proof}

Another important observation is that the Frobenius commutes with
\'etale pull back, i.e. $F_*f^*\cal{M} \cong f^*F_*\cal{M}$ for all
$\OO_X$-modules $\cal{M}$.

\begin{lemma}
Let $X$ and $Y$ be $k$-varieties and $f: Y \to{} X$ an \'{e}tale morphism. Then
\[
\xymatrix{ Y \ar^F[r] \ar_f[d] & {Y'} \ar^f[d] \\
            X \ar^F[r]          & {X'}
            }
\]
is a fiber square, i.e. $Y \cong {Y'} \times_{X'} X$.
\end{lemma}
\begin{proof} First assume that $f$ is finite and \'etale. Then by the
above diagram and the universal property of the fiber product we have
a natural map $Y \to{} Y' \times_{X'} X $. To check this map is an
isomorphism one reduces to the case where $X$ (and therefore $Y$) is
affine. Since $Y' \times_{X'} X \cong Y \times_{F_X} X$ one further
reduces to the case of the absolute Frobenius. Thus it remains to show
that for a finite \'etale morphism of $k$-algebras $R \to{} S$ we have $(F_S)_* S
\cong S \otimes_{F_R} (F_R)_*R$. This is done in \cite[page 50]{HH90}.
\\ In general an \'etale map can be factored into a finite \`etale map
followed by an open immersion: $f:Y \to{} U \to{} X$. It remains to
check the statement for the open immersion $U \to{} X$. As before one
reduces to the case of the absolute Frobenius. As topological spaces
$U \times_{F_X} X$ and $U$ are the same and a simple local calculation
shows that they are isomorphic as $k$-schemes.
\end{proof}

With this at hand we can apply \cite[Proposition III.9.3]{Hartshorne} and get:

\begin{corollary}\label{etaleFrob}
Let $\cal{M}$ be an $\OO_X$-module and $f: Y \to{} X$ \'etale. Then $f^*F_*
\cal{M} \cong F_*f^* \cal{M}$.
\end{corollary}

Since similar in spirit to the above and used later we note here a Lemma which
shows that $F_*$ commutes with completion. For this we consider the local case.

\begin{lemma}\label{completeFrob}
Let $(R,m)$ be a local $k$-algebra, essentially of finite type over $k$. Then
$\hat{R} \cong R \otimes_{R'} \hat{R}'$. Here, $\hat{R}$ denotes the $m$-adic
completion of $R$.
\end{lemma}
\begin{proof}
As in the previous Lemma $R \otimes_{R'} \hat{R}' \cong R \otimes_{F_R}
\hat{R}$ and we reduce to the case of the absolute Frobenius. Now it follows
from the definition of completion:
\[
R \otimes_{F_R} \hat{R} = R \otimes_{F_R} \varprojlim \frac{R}{m^i} \cong
\varprojlim \frac{R}{R \otimes_{F_R}m^i} = \varprojlim \frac{R}{(m^i)^{[p]}R} =
\hat{R}
\]
In the last equality we used that the sequence $(m^i)^{[p]}$ is cofinal with
$m^i$ and therefore they both have the same limit.
\end{proof}

An easy calculation with the equation of the last Lemma shows that for any
$R$-module $M$ $F_*$ commutes with completion, i.e.
$F_*\hat{M} \cong \hat{R}' \otimes_{R'} F_*M$.

\subsection{Cartier operator on Zariski differentials}

The Cartier isomorphism on the smooth scheme $U'$ gives after applying $i_*$ an
isomorphism
\[
   i_*C^{-1}:\ \Omtilde^a_{X'} \xrightarrow{\
\cong\ \ } i_* \HH^a(F_*\Omega^{\bullet}_U).
\]
To obtain an analog of the Cartier operator for the Zariski sheaves it remains
to construct natural maps
\[
  \phi_a:\ \HH^a(F_*\Omtilde^{\bullet}_X) \To i_* \HH^a(F_*\Omega^{\bullet}_U)
\]
for all $a \geq 0$. This is obtained with the help of the following general
Lemma.

\begin{lemma}
Let $\cal{A}^\bullet$ be a complex in an abelian category and $G$ a left exact
functor of abelian categories. One has natural maps
\[
  \phi_a: \ \HH^a(G\cal{A}^\bullet) \to{} G\HH^a(\cal{A}^\bullet)
\]
for all $a \geq 0$.
\end{lemma}
\begin{proof} To see this denote by $\cal{Z}^a$, $\cal{B}^a$ the cycles and
boundaries of $\cal{A}^\bullet$ (resp. $\cal{Z}^a_G$, $\cal{B}^a_G$ for
$G\cal{A}^\bullet$). By the left exactness of $G$ we get $\cal{Z}^a_G \cong
G\cal{Z}^a$. Applying $G$ to the sequence $0 \to{} \cal{Z}^a \to{} \cal{A}^a
\to{} \cal{B}^{a+1} \to{} 0$ and arranging with the equivalent sequence for
$G\cal{A}^\bullet$ in the diagram
\[
  \xymatrix{
  {0} \ar[r] &{\cal{Z}^a_G} \ar@{=}[d] \ar[r] &{G\cal{A}^a} \ar@{=}[d] \ar[r] &{\cal{B}^{a+1}_G} \ar@{_(-->}[d]\ar[r] & {0} \\
  {0} \ar[r] &{G\cal{Z}^a} \ar[r] &{G\cal{A}^a} \ar[r] &{G\cal{B}^{a+1}}
  }
\]
we see that (Snake Lemma) the vertical arrow to the right is an injection.
Now, using the equivalent diagram for the short exact sequence defining
cohomology
\[
  \xymatrix{
  {0} \ar[r] &{\cal{B}^a_G} \ar@{_(->}[d] \ar[r] &{\cal{Z}^a_G} \ar@{=}[d] \ar[r] &{\HH^a(G\cal{A}^\bullet)} \ar@{-->}[d]^{\phi_a} \ar[r] &{0} \\
  {0} \ar[r] &{G\cal{B}^a} \ar[r] &{G\cal{Z}^a} \ar[r] &{G\HH^a(\cal{A}^\bullet)}
  }
\]
together with the Snake Lemma gives the dotted vertical arrow to the right.
Furthermore, $\phi_a$ is injective if and only if $\cal{B}^a_G \to{}
G\cal{B}^a$ is an isomorphism and surjective if and only if the map $G\cal{Z}^a
\to{} G\HH^a(\cal{A}^\bullet)$ is surjective.
\end{proof}

Since $i_*$ is left exact we can apply this construction to the complex
$F_*\Omega^\bullet_U$ and therefore obtain natural maps $\phi_a :
\HH^a(i_*F_*\Omega^{\bullet}_U) \To i_* \HH^a(F_*\Omega^{\bullet}_U)$ for all
$a \geq 0$. Since $F_*$ commutes with $i_*$ the first term is just
$\HH^a(F_*\Omtilde^{\bullet}_X)$

By the Cartier isomorphism on the smooth locus $\HH^a(F_*\Omega_U^\bullet)$ is
isomorphic to $\Omega^a_{U'}$ which is a locally free and therefore reflexive
$\OO_{U'}$-module. By Lemma \ref{ReflExt} $i_*\HH^a(F_*\Omega_U^\bullet)$ is a
reflexive $\OO_{X'}$-module. Thus, by the universal property of the reflexive
hull, the map $\phi_a$ factors through the natural map from
$\HH^a(F_*\Omtilde_X^\bullet)$ to its double dual:
\[
  \xymatrix{
  {\HH^a(F_*\Omtilde_X^\bullet)} \ar^{\phi_a}[r] \ar[d]_{(\ )^{**}}
  &{i_*\HH^a(F_*\Omega_U^\bullet)} \\
  {\HH^a(F_*\Omtilde_X^\bullet)^{**}} \ar_{\cong}[ru]
  }
\]
By the uniqueness part of Lemma \ref{ReflExt} the diagonal map has to be an
isomorphism since $\phi_a$ is an isomorphism on $U'$. If we use this natural
isomorphism to identify $i_*\HH^a(F_*\Omega_U^\bullet)$ with the double dual of
$\HH^a(F_*\Omtilde_X^\bullet)$ the map $\phi_a$ becomes nothing but double
dualizing. Composing $\phi_a$ with the inverse of $i_*C^{-1}$
 we get:

\begin{proposition}
  For all $a \geq 0$ the Cartier isomorphism on $U$ induces a natural
  $\OO_{X'}$-linear map
  \[ \wt{C}:\ \HH^a(F_*\Omtilde^{\bullet}_X) \To \Omtilde^a_{X'}. \]
\end{proposition}

By construction, $\wt{C}$ is an isomorphism if and only if the map $\phi_a$ is.
This in turn is the case if and only if $\HH^a(F_*\Omtilde^\bullet_X)$ is
reflexive. This observation is key for the proof of the next Proposition.

\begin{proposition}\label{etaleCartier}
  Let $f: Y \to{} X$ \'etale. If the Cartier operator $\wt{C}_X$ on $X$ is an
  isomorphism then so is $\wt{C}_Y$ on $Y$. If in addition $f$ is
  finite, then the reverse implication also holds.
\end{proposition}
\begin{proof}
With the above observation we only have to show that if
$\HH^a(F_*\Omtilde^\bullet_X)$ is reflexive then so is
$\HH^a(F_*\Omtilde^\bullet_Y)$. Since an \'etale map is flat we can
apply Lemma \ref{refl:flat} to see that the reflexivity of
$\HH^a(F_*\Omtilde^\bullet_X)$ implies that
$f^*\HH^a(F_*\Omtilde^\bullet_X)$ is reflexive. By flatness $f^*$ commutes with
taking cohomology and thus
\[
  f^*\HH^a(F_*\Omtilde^\bullet_X)=\HH^a(f^*F_*\Omtilde^\bullet_X).
\]
By Corollary \ref{etaleFrob} $f^*F_* \Omtilde^a_X = F_*f^*\Omtilde^a_X$ which
by \ref{Zariskietale} is just $F_* \Omtilde^a_Y$. \\
If $f$ is also finite then it is faithfully flat and an application of
Corollary \ref{refl:fflat} yields the reverse implication.
\end{proof}

\section{Cartier isomorphism for toric varieties}

In this section we give the proof that the Cartier map is an isomorphism for
affine toric varieties over $k$. The main ingredient is a description of the
Zariski-\deRham complex due to Danilov \cite{Danilov78}. An interesting aspect
of the given proof is that it is analogous to Danilov's proof of a Poincar\'e
type lemma for toric varieties which says that the Zariski-\deRham complex is a
resolution of the constant sheaf $k$. Thus, as in the smooth case, the Cartier
isomorphism can be viewed as a characteristic $p$ analog of the Poincar\'e
lemma, both determining the cohomology of the \deRham complex. Furthermore,
using the results of the previous section on the stability of the Cartier
isomorphism under \'etale morphisms we conclude that also for toroidal
varieties the Cartier operator is an isomorphism.

\subsection{Affine toric varieties}

First we recall some notation from toric geometry, for a more detailed
introduction to toric geometry see \cite{FultonToric} or \cite{Danilov78}. Let
$N \cong \ZZ^n$ be a lattice, $M=N^*$ the dual lattice and we denote $M
\otimes_{\ZZ} \QQ$ by $M_\QQ$. The identification of $N$ with $\ZZ^n$ allows us
to consider the elements of $N$ as $n$-tuples of integers, similarly for $M$.

By a cone in $N$ we mean a subset $\sigma \subseteq N_{\QQ}$ of the form
$\sigma=\{r_1v_1+ \dotsb +r_sv_s|r_i \geq 0\}$ for some $v_i$ in $N$. The dual
cone of $\sigma$ in $M$ is $\sigma^{\vee}=\{u \in M_{\QQ}| (u,v) \geq 0, \text{
for all } v\in \sigma\}$ where $(\ ,\ )$ is induced from the natural
pairing between the dual lattices $M$ and $N$. A face of $\sigma$ is any set
$\sigma \cap u^{\bot}$ for some $u \in \sigma^{\vee}$ and is a cone in $N$. For
cones $\tau$ and $\tau'$ in $M_{\QQ}$ $(\tau-\tau')$ denotes the set of all
$t-t'$ for $t \in \tau$ and $t' \in \tau'$.

Let $\sigma^{\vee}$ be a cone in $M_{\QQ}$. Let $k$ be a perfect field and let
$A=k[\sigma^{\vee} \cap M]$ be the affine semigroup ring corresponding to
$\sigma$, i.e. $A$ is the $k$-subalgebra of $k[x_1,x^{-1}_1, \dots
,x_n,x_n^{-1}]$ generated by the monomials $x^m=x_1^{m_1}\cdot \ldots \cdot
x_n^{m_n}$ for $m \in \sigma^{\vee}$. $X=X_\sigma = \Spec(A)$ denotes the
corresponding affine toric variety. Note that $A$ carries a natural $M$-grading
where the generator $x^m$ corresponding to $m \in \sigma^{\vee}
\subseteq M$ of $A$ is given degree $m$. As is well known $X$ is
normal \cite{FultonToric} and we have the notion of Zariski
differentials defined above. As before we denote by $U$ its smooth
locus. Since $X$ is affine $\Omtilde^a_X$ is the sheaf
associated to some $A$-module, which we will call $\Omtilde^a_A$.

\subsubsection{Danilov's description}

There is a nice description of the modules $\Omtilde^a_A$ due to Danilov
\cite{Danilov78} which we will give next.

Let $V$ denote the vector space $M \otimes_{\ZZ} k$. For each codimension one
face $\tau$ of $\sigma^{\vee}$ define a subspace $V_{\tau} \subset V$.
\[
  V_{\tau}=(M \cap (\tau-\tau))\otimes_{\ZZ}k
\]
For every $m \in \sigma^{\vee}$ we define a subspace $V_m \subseteq V$:
\[ V_m = \bigcap_{\tau \ni m} V_{\tau} \]
where $\tau$ is ranging over the codimension one faces of $\sigma^{\vee}$
containing $m$.

\begin{proposition}\label{omegaA}
  $\Omtilde^a_A$ is isomorphic to the graded $A$-submodule $V^a_\sigma$ of the $M$-graded
  $A$-module $\bigwedge^a V \otimes_k A$ with degree $m \in M$ piece
  \begin{equation*}
    V^a_\sigma(m)=\bigwedge^a V_m \cdot x^m .
  \end{equation*}
  The isomorphism is given by the following degree preserving map of graded
  $A$-modules
  \[ \xymatrix@R=0pt{
    {{f_a:\ \ \ V_\sigma^a}} \ar[r]  &{\Omtilde^a_A} \\
    {\qquad \qquad \scriptstyle (m_1\otimes1\wedge \dots \wedge m_a\otimes 1) \cdot x^m}
    \ar@{|->}[r] &{\scriptstyle x^{m-m_1\dots-m_a}dx^{m_1}\wedge \dots \wedge
    dx^{m_a}}}
  \]
\end{proposition}

For the proof see \cite[p. 110]{Danilov78}. With this identification of
$\Omtilde_A^a$ with $V_\sigma^a$ the differential $d:\Omtilde^a_A \to{}
\Omtilde^{i+1}_A$ on the piece of degree $m$ is given by wedging with $m
\otimes 1$. Thus the degree $m$ part of the complex $\Omtilde^{\bullet}_X$ is
just
\begin{equation*}
\xymatrix@C=.5cm{
  0 \ar[r] & {\bigwedge^0 V_m} \ar[rr]^{m\otimes 1 \cdot} && {\bigwedge^1 V_m} \ar[rr]^{m\otimes 1 \wedge} && {\bigwedge^2 V_m} \ar[rr]^{m\otimes 1 \wedge} && {\hspace{12pt}\cdots\hspace{12pt}}}
\end{equation*}
The key observation is that this complex is exact whenever $m \otimes 1$ is a
nonzero element of the vector space $V_m$.

If the characteristic of $k$ is zero this is the case if and only if
$m \neq 0$.
Thus the only degree on which the complex is not exact is degree $0 \in M$.
Assuming that the cone $\sigma^\vee$ does not contain any linear subspace we
see that $V_0=0$ and therefore in degree zero we just get $\bigwedge^0 V_0=k$.
This was Danilovs proof of the Poincar\'e lemma for the Zariski-\deRham complex
on an affine toric variety:

\begin{PL}[\cite{Danilov78}] If the characteristic of $k$ is
zero and $\sigma^{\vee}$ is a cone with vertex (i.e. $\sigma^{\vee}$ does not
contain a linear subspace) then
\[ \xymatrix@C=0.5cm{ 0 \ar[r] & {k} \ar[r]^{} & {\Omtilde^{\bullet}_A}} \]
is exact, i.e. $\Omtilde^{\bullet}_A$ is a resolution of $k$.
\end{PL}

Back to the case when the characteristic of $k$ is greater than zero. Then we
get the following analog of \ref{smooth:cartier} for affine toric varieties.

\begin{theorem}\label{toriccartier} If $k$ has characteristic $p > 0$ then
the map
\[
  \phi: \bigoplus \Omtilde^{a}_{A'}[-a] \To F_*\Omtilde^{\bullet}_A
\]
sending $\ (m'\otimes 1) \cdot x^m\ $ to $\ (m'\otimes 1) \cdot x^{pm}\ $ is a
split injection and induces an isomorphism $C^{-1}: \Omtilde^a_{A'} \to{}
\HH^a(F_*\Omtilde^{\bullet}_A)$.
\end{theorem}
\begin{proof} As already in the statement of this Theorem we make
use of the identification in \ref{omegaA} without further mentioning.

By definition $\phi$ is additive. Noting that $A' = A \otimes_{F_k} k$ and
$\Omtilde^a_{A'} = \Omtilde^a_A \otimes_{F_k} k$ it is easily checked that
$\phi$ is also $A'$-linear. Since $m \otimes 1 = 0$ if and only if $m \in pM$
we see that the cohomology lives exclusively in degrees $pm$ for $m \in M$; in
these degrees the differential (wedging with $pm \otimes 1 = 0$) is the zero
map. $\phi$ maps the graded piece of degree $m$ isomorphically to the graded
piece of degree $pm$. Thus the image of $\phi$ is exactly the graded pieces of
degree $pm$ for $m \in M$ which is the cohomology of $F_*\Omtilde^{\bullet}_A$.
Thus $\HH^a(F_*\Omtilde^{\bullet}_A)$ is isomorphic to the graded
$A'$-submodule of degree $pM$ pieces of $F_*\Omtilde^a_A$ and therefore $\phi$
is also split. Using \ref{omegaA} one easily verifies that $d \circ \phi=0$,
i.e $\phi$ is a map of complexes.

It remains to show that the above defined map actually is an inverse of the
Cartier operator. Using the description of \ref{omegaA} we see that for all $m
\in \sigma^{\vee}$ $\ \phi(dx^m \otimes 1)=x^{(p-1)m}dx^m$ which is just the
defining property of the inverse Cartier map from which $\wt{C}$ was
constructed.
\end{proof}

Since checking that the Cartier operator $\wt{C}$ is an isomorphism is a local
matter we get as an immediate Corollary:

\begin{corollary}
Let $X$ be a toric variety over a perfect field $k$. Then the Cartier operator
\[
  \wt{C}:\ \HH^a(F_*\Omtilde^{\bullet}_X) \To \Omtilde^a_{X'}
\]
is an isomorphism.
\end{corollary}

Theorem \ref{toriccartier} for \emph{affine} toric varieties is much stronger
than what we just established for the toric case. Not only is the Cartier map
an isomorphism but also it is induced by a split injection $\Omtilde^a_{X'}
\to{} F_* \Omtilde^a_X$. This stronger statement can still be achieved for
toric varieties in general (not necessarily affine). The key point is that one
has to ensure that the map $\phi$ defined locally on the affine toric pieces as
above patches to give a global map of sheaves.

For simplicity of notation we assume in the following discussion that
$k=\FF_p$. Notice that the description of $\Omtilde^a_X$ in Lemma \ref{omegaA}
localizes. I.e. if $\tau$ is a face of $\sigma$ then $X_\tau$ is an open
subset of $X_\sigma$ and we have an inclusion $\Omtilde^a_{X_\sigma} \subseteq
\Omtilde^a_{X_\tau}$. This corresponds to the inclusion $V^a_\sigma \subseteq
V^a_\tau$. The reason for this is that the map $V^a_\sigma \To
\Omtilde^a_{X_\sigma}$ is induced from the map $\bigwedge^a V \to{}
\Omtilde^a_{T_N}$ for \emph{any} cone $\sigma$ in $N$, here $T_N$ is the
$n$-dimensional torus.

Now let $X$ be the toric variety associated to a fan $\Sigma$ in $N$ (cf.
\cite{FultonToric}). The torus $T_N$ is an open subset of $X$ and also of each
affine toric piece $X_\sigma$ of $X$ for $\sigma \in \Sigma$. Thus
$\Omtilde^a_{X_\sigma} \subseteq \Omtilde^a_{T_N}$ and the map
$\phi_\sigma:\Omtilde^a_{X_\sigma} \to{} F_*\Omtilde^a_{X_\sigma}$ is induced
from the corresponding map for the torus $T_N$ ($T_N$ is the toric variety
associated the cone $\{0\} \subseteq N$). This shows that the maps
$\phi_\sigma$ agree on the intersections of their domains and thus give rise to
a globally defined map $\phi$ as required. We just proved (well, sketched a
proof of) the following Proposition.

\begin{proposition}
Let $X$ be a toric variety over a perfect field $k$. Then the Cartier
operator is an isomorphism and is induced by a split injection
\[
  0 \to{} \bigoplus_a \Omtilde^a_{X'}[-a] \to{} F_*\Omtilde^\bullet_X.
\]
\end{proposition}

The existence of a split injection inducing the Cartier operator was
already shown by Buch, Thomsen, Lauritzen and Mehta  in \cite{Anders}
using the existence of natural liftings
of characteristic $p$ toric varieties to toric varieties over
$\ZZ/p^2\ZZ$, along with liftings of Frobenius. Their methods though do
not yield any information on whether the Cartier operator is an
isomorphism, but are sufficient to imply Bott vanishing and
the degeneration of the Hodge to \deRham spectral sequence for projective toric
varieties.

\subsection{Toroidal varieties}
With the help of Proposition \ref{etaleCartier} we can extend the results from
the previous section to a important class of varieties: toroidal embeddings.

\begin{definition}
A normal scheme $X$ locally of finite type over a field $k$ will be called
\emph{weakly toroidal} if for each point $x \in X$ there is a neighborhood
$U(x)$ of $x$ and a diagram
\[
\xymatrix{       & {Y(x)} \ar_{f}[ld] \ar^g[rd] &  \\
          {Z(x)} &        & {U(x)}}
\]
where $Z(x)$ is an affine toric variety, $f$ is \'etale and $g$ is
finite and \'etale.
\end{definition}

This class of varieties includes the class of toroidal embeddings as defined in
\cite{TE}. We have the following Theorem.

\begin{theorem}
Let $X$ be weakly toroidal. Then the Cartier operator $\wt{C}$ on $X$ is an
isomorphism.
\end{theorem}
\begin{proof} To see that the Cartier operator on $X$ is an isomorphism is a local
issue so we can assume that there is a diagram
\[
\xymatrix{       & {Y} \ar_{f}[ld] \ar^g[rd] &  \\
          {Z} &        & {X}}
\]
with $Z$ affine toric, both maps \'etale and $g$ also finite. On the
toric variety $Z$ the Cartier operator is an isomorphism and thus by
\ref{etaleCartier} it is an isomorphism on $Y$, and again by
\ref{etaleCartier}, it is an isomorphism on $X$.
\end{proof}

\section{An Obstruction to surjectivity}

The last section shows that the Cartier operator is an isomorphism in the case
of toroidal embeddings and it is natural to ask what is the exact class of
varieties for which the Cartier operator is an isomorphism. It seems too much
to hope that it would be an isomorphism for all normal varieties. And in fact
there is a necessary and sufficient condition for $\wt{C}$ to be a surjection
at the top level, i.e. for $\wt{C}: \HH^d(F_*\Omtilde^{\bullet}_X) \to{}
\Omtilde^d_X$ where $d$ is the dimension of $X$.

To simplify the notation we assume $k = \FF_p$. Let $(R,m)$ be a normal
Cohen-Macaulay local domain of dimension $d$, essentially of finite type over
$k$. Then $R$ is called $F$-injective if the Frobenius $F: R \to{} F_*R$
induces an injection on local cohomology $H^d_m(R) \to{} H^d_m(F_*R)$. To
relate this to the Cartier operator on $X=\Spec R$ we denote the global
sections of $\Omtilde^d_X$ by $\omega_R$ which by normality is a canonical
module for $R$. Since the Cartier operator $\Gamma(X,\HH^d(F_*\Omtilde^\bullet_X)) \to{}
\omega_R$ is surjective if and only if its composition with the natural
surjection $F_*\omega_R \to{} \Gamma(X,\HH^d(F_*\Omtilde^\bullet_X)) \to{} \omega_R$ is
surjective we consider now the map $C: F_*\omega_R \to{} \omega_R$ to which we
also refer to as the Cartier map. This map is surjective iff its $m$-adic
completion $\hat{C}: \hat{R} \otimes F_*\omega_R \to{} \hat{R} \otimes
\omega_R$ is a surjective map of $\hat{R}$-modules by the faithfully flatness
of completion. By the remark after Lemma \ref{completeFrob} $F_*$ commutes with
completion and by \cite[3.3.5]{BrunsHerzog} $\ \hat{R} \otimes \omega_R$ is a
canonical module for $\hat{R}$. Thus $\hat{R} \otimes \omega_R \cong
\omega_{\hat{R}}$ and $\hat{R} \otimes F_*\omega_R \cong F_*\omega_{\hat{R}}$.
Let $E$ denote the injective hull of the residue field of $\hat{R}$ and take
the Matlis dual of the completed Cartier map $\hat{C}$ to get
\[
\Hom_{\hat{R}}(\omega_{\hat{R}},E) \to{} \Hom_{\hat{R}}(F_*\omega_{\hat{R}},E).
\]
Note that, since $F$ is a finite map $F_*\omega_{\hat{R}} \cong
\Hom_{\hat{R}}(F_*\hat{R},\omega_{\hat{R}})$ (cf. \cite[3.3.7]{BrunsHerzog}).
This together with the local duality theorem \cite[3.5.8]{BrunsHerzog} shows
that the last map is just
\[
  H^d_{\hat{m}}(\hat{R}) \to{} H^d_{\hat{m}}(F_*\hat{R})
\]
where we can leave out the completion and get the map induced by the Frobenius
$F: H^d_m(R) \to{} H^d_m(F_*R)$. Summarizing we get the following Proposition.
\begin{proposition}
Let $(R,m)$ be a local Cohen-Macaulay domain essentially of finite type over $k$. Then the
Cartier map $C: F_*\omega_R \to{} \omega_R$ is surjective if and only if $R$ is
$F$-injective.
\end{proposition}
One implication of this Proposition was already noticed by Mehta and
Shrinivas in \cite{MehtaSr}  where they use the Cartier operator on
the top spot to obtain information about the singularities of the variety.

\bibliographystyle{alpha}

\end{document}